\documentclass[a4paper, 10pt, reqno]{amsart}

\usepackage{ amssymb, amsmath, amsthm}
\usepackage{graphicx, psfrag, setspace, subfigure, gensymb}
\usepackage{amsfonts}
\usepackage{bbm}
\usepackage{fix-cm}
\usepackage{comment}
\usepackage{tikz,tikz-cd}
\usetikzlibrary{matrix,arrows,decorations.pathmorphing,decorations.pathreplacing}
\tikzset{commutative diagrams/diagrams={baseline=-2.5pt},commutative diagrams/arrow style=tikz}
\usepackage[colorlinks]{hyperref}
\usepackage{float}
\usepackage{setspace} 

\newcommand\Z{\mathbb Z}
\newcommand\C{\mathbb C}

\newcommand{\cC}{\mathcal{C}}
\newcommand{\cD}{\mathcal{D}}
\newcommand{\cE}{\mathcal{E}}

\newcommand{\cL}{\mathcal{L}}

\newcommand{\cO}{\mathcal{O}}

\newcommand\id{\mathrm 1}

\newcommand\into{\hookrightarrow}
\newcommand\To{\longrightarrow}

\newcommand\Hom{\operatorname{Hom}}

\renewcommand\hom{\mathcal{H}om }

\newcommand\Ext{\operatorname{Ext}}
\renewcommand\Im{\operatorname{Im}}

\renewcommand\P{\mathbb P}

\newcommand\Perf{\operatorname{Perf}}

\newcommand{\al}[1]{\begin{align*}#1\end{align*}}

\newcommand{\beq}[1]{\begin{equation}\label{#1} }
\newcommand{\eeq}{\end{equation}}
\newcommand{\pgap}{\vspace{5pt}}

\newtheorem{prop}[equation]{Proposition}
\newtheorem{thm}[equation]{Theorem}

\theoremstyle{remark}
\newtheorem{rem}[equation]{Remark}
\newtheorem*{acks}{Acknowledgements}
\newtheorem*{conv}{Conventions}

\theoremstyle{definition}
\newtheorem{defn}[equation]{Definition}
\newtheorem{eg}[equation]{Example}

\makeatletter \@addtoreset{equation}{section} \makeatother

\let\oldtocsection=\tocsection
\let\oldtocsubsection=\tocsubsection
\let\oldtocsubsubsection=\tocsubsubsection
\renewcommand{\tocsection}[3]{\hspace{0em}\oldtocsection{#1}{#2}{#3}}
\renewcommand{\tocsubsection}[3]{ \hspace{1em} \oldtocsubsection{#1}{\small{#2}}{\small{#3}} }
\renewcommand{\tocsubsubsection}[3]{\hspace{2em}\oldtocsubsubsection{#1}{\small{#2}}{\small{#3}}}

\newcommand{\marginparstretch}{0.6}
\let\oldmarginpar\marginpar
\renewcommand\marginpar[1]{\-\oldmarginpar[\framebox{\setstretch{\marginparstretch}\begin{minipage}{\marginparwidth}{\raggedleft\scriptsize #1}\end{minipage}}]{\framebox{\setstretch{\marginparstretch}\begin{minipage}{\marginparwidth}{\raggedright\scriptsize #1}\end{minipage}}}}

\newcommand{\qquotes}[1]{``{#1}''}

\newcommand{\aand}{\quad\quad\mbox{and}\quad\quad}

\newcommand{\D}{D^b}
\newtheorem{'prop'}[equation]{\qquotes{Proposition}}

\begin{document}

\title{All autoequivalences are spherical twists}%
\author{Ed Segal}%

\maketitle

\begin{abstract}
In this short note we observe that, for purely formal reasons, any autoequivalence can be constructed as a twist around a spherical functor.  As an example, we show how the P-twists constructed by Huybrechts and Thomas can be formulated as spherical twists.
\end{abstract}

\tableofcontents

\section{Introduction}

In \cite{ST}, Seidel and Thomas introduced a new symmetry of derived categories called spherical twists. These are formed by a simple construction from a \emph{spherical object}, which is an object whose self-Ext algebra looks like the cohomology of a sphere.  Since then various authors (Horja, Rouquier, Toda, Anno, Logvinenko,...) have generalized this construction by replacing a spherical object with a \emph{spherical functor}, which is a functor
$$F: \cC \to \cD $$
satisfying certain axioms. We refer the reader to \cite{Add} for a survey of the development of this idea (including references to the authors named above), and to \cite{AnnLog} for a very general technical treatment. Roughly, we should think of a spherical functor as a family of objects in $\cD$ parametrized by the source category $\cC$, such that each object is spherical relative to the base. If the source category $\cC$ is $\D(pt)$ then $F$ is exactly a choice of spherical object in $\cD$.

 If $R$ is the right adjoint of $F$ then by forming the cone on the counit we get an endofunctor of $\cD$
$$T = [FR \to \id_\cD]$$
and the important fact is that $T$ is an autoequivalence. This is called the \emph{spherical twist} around $F$. 

 The purpose of this short note is to observe that in fact \emph{any} autoequivalence of a derived category can be formulated as a spherical twist. Given the data of a category $\cD$ and an autoequivalence $\Phi$, we can formally construct a source category which we call $\cD_\Phi$  and a spherical functor which we call
$$j_*: \cD_\Phi \to \cD$$
 such that $\Phi$ is the  spherical twist around $j_*$. This construction is explained in Section \ref{sect.general}.

If our autoequivalence is given to us as a spherical twist around some spherical functor $F:\cC \to \cD$ then the construction of Section \ref{sect.general} does not recover $\cC$ and $F$. Section \ref{sect.withsection} gives a sketch proof of a solution to this issue. We argue that (under a necessary hypothesis) the pair $(\cC, F)$ is in fact a deformation of $(\cD_\Phi, j_*)$. This section is only a sketch because the construction suffers from foundational issues that we have not managed to resolve.
\pgap

 In some sense the result of this paper is a negative one - it says that there is nothing special about spherical twists in the abstract. However, in practice they remain a useful way of constructing complicated autoequivalences (of $\cD$) out of simpler ones (of $\cC$). 

The result also suggests that if we encounter a new kind of autoequivalence then it may be useful to describe it as a spherical twist, particularly if we can find one with a `minimal' source category (the spherical functor constructed in Section \ref{sect.general} is so tautological that it seems unlikely to be enlightening). In Section \ref{sect.Ptwists} we illustrate this philosophy by showing that the $\P$-twists defined by Huybrechts and Thomas \cite{HuyTho} can be formulated as spherical twists, with a very simple source category.

\begin{acks} I would like to thank Will Donovan, Tim Logvinenko, and especially Michael Groechenig for many helpful conversations. Nick Addington suggested some major improvements which ended up delaying the completion of the paper for a long time, for this I would like to express my gratitude. I guess.
\pgap

This project has received funding from the European Research Council (ERC) under the European Union Horizon 2020 research and innovation programme (grant agreement No.725010).
\end{acks}

\begin{conv}
We work with triangulated categories and dg-categories defined over some field $k$. All functors between triangulated categories are assumed to be exact and $k$-linear, and when we write $\Hom$, $\otimes$, etc. we mean the derived versions of these functors.
\end{conv}

\section{General autoequivalences}\label{sect.general}

Let $\cC$ and $\cD$ be triangulated categories, and let 
$$F: \cC \to \cD$$
be a functor with right and left adjoints $R$ and $L$. 
\begin{defn}\label{defn:spherical}We say that $F$ is \emph{spherical} if any two of the following conditions hold.
\begin{list}{(\roman{enumi})}{\usecounter{enumi}}
\item The \emph{twist} functor
$$T := \big[FR \to \id_\cD\big]$$
(the cone on the counit of the adjunction) is an autoequivalence of $\cD$.
\item The \emph{cotwist} functor
$$ C := \big[\id_\cC \to RF\big][-1]$$
(the shifted cone on the unit of the adjunction) is an autoequivalence of $\cC$.
\item $R$ and $L$ are related by:
$$ R = LT[-1]$$
\item $R$ and $L$ are related by:
$$ R = CL[1] $$
\end{list}
\end{defn}
 The utility of this definition is that if any two of these conditions hold, then the other two hold automatically \cite{AnnLog} (see also \cite{Meachan}). When $F$ is spherical we call $T$ a spherical twist.

\begin{rem}\label{rem.foundations1} We are being a little careless here: to be able to take cones on functors we should really choose dg-enhancements of $\cC$ and $\cD$ and assume we have lifts of $F, R$ and $L$ to quasi-functors (see for example \cite[\S 2]{AnnLog} or \cite[\S 4]{Toen}). Alternatively we can assume we're in some geometric or algebraic context where we have a sensible model for the functor category (Fourier-Mukai kernels, bimodules, etc).
\end{rem}

 We wish to show that any autoequivalence of $\cD$ arises from some spherical functor. We will begin by considering two (closely analogous) examples, and extrapolate from there to the general case.

\subsection{Two examples}\label{sect.2examples}

\begin{eg}\label{eg.algebraic} Let $A$ be an algebra, commutative or non-commutative, and let $\cD = \D(A)$ be the derived category of right $A$-modules. Let $B$ be an invertible bimodule over $A$, so the endofunctor
$$\otimes B: \D(A) \to \D(A) $$
is an autoequivalence. Let us show that $\otimes B$ can be formulated as a spherical twist.

Let $B^\vee = \Hom_{-A}(B,A)$ denote the inverse bimodule, and form the extension algebra
$$E = A\oplus B^\vee$$
 (by declaring that the product of two elements in $B^\vee$ is zero). This comes with an algebra map $j: A \to E$, and an associated functor:
$$j_*: \D(E) \to \D(A) $$
The left adjoint of $j_*$ is $j^* = \otimes_A E$, and the right adjoint of $j_*$ is: 
$$j^! = \Hom_A(E, -)  = \otimes_A(A\oplus B) = j^*\circ(\otimes B)$$
The twist around $j_*$ is the functor:
$$T = \otimes B [1] $$
Conditions (i) and (iii) from Definition \ref{defn:spherical} are satisfied, so $j_*$ is spherical and $T$ is a spherical twist.

 We have picked up an unwanted shift, but this is easily corrected. Instead of $A\oplus B^\vee$ we should form the graded algebra $A\oplus \left( B^\vee[1]\right)$, and view it as a dga with no differential. Then $j$ becomes a map of dga's, and the twist around $j_*$ is precisely $\otimes B$.

\end{eg}

\begin{eg}\label{eg.geometric} Let $X$ be a scheme, and let $\cD = \D(X)$, the derived category of coherent sheaves on $X$. Now let $\cL$ be a line-bundle on $X$, so we have an autoequivalence:
 $$\otimes \cL: \D(X) \to \D(X)$$
Up to a shift, we will formulate this autoequivalence as a spherical twist.  

We consider the total space of the dual line-bundle $\cL^{-1}$, with its zero section: 
$$i: X \into\cL^{-1}$$
 We consider the functor:
$$ i^*: \D(\cL^{-1}) \to \D(X)$$
The right adjoint of $i^*$ is $i_*$, and the left adjoint is $(\otimes \cL^{-1})\circ i_*[-1]$. Then the twist around $i^*$ is the functor:
$$T = \big[ i^*i_* \to \id_{\D(X)} \big]$$
For any object $\cE\in \D(X)$ we have that 
$$i^*i_*\cE = \cE \oplus (\cE\otimes \cL[1])$$
and hence:
$$T =  \otimes\cL[2] $$
Again conditions (i) and (iii) from Definition \ref{defn:spherical} are satisfied, so the functor $i^*$ is spherical and $T$ is a spherical twist.

Notice that we don't really need the whole category $\D(\cL^{-1})$ here, the construction only uses the essential image of the functor $i_*$. This is the subcategory $\D_X(\cL^{-1})\subset \D(\cL^{-1})$  consisting of objects which are set-theoretically supported on the zero section. A similar remark applies in the previous example.
\end{eg}

The intersection of these two examples is the case where $A=\cO_X$ is a commutative algebra and $B=\cL$ is an invertible bimodule for which the left and right actions coincide - this is a line bundle on an affine scheme. To see the connection between the two constructions we apply Koszul duality in Example \ref{eg.geometric}, as follows. Let's use  $\D_X(\cL^{-1})$ as our source category.  This is equivalent to the derived category of
$$ \cE:= i^*\hom_{\cL^{-1}}(i_*\cO_X, i_*\cO_X) = \cO_X \oplus \left( \cL^{-1}[-1]\right)$$
which is  a sheaf of graded algebras over $X$, \emph{i.e.}~the $\cO_X$-linear functor
$$\otimes_\cE  i_*\cO_X: D^b(mod-\cE) \to D^b(\cL^{-1})$$
is an embedding, with image $D^b_X(\cL^{-1})$. In an affine neighbourhood in $X$ this is clear, since then $i_*\cO_X$ generates the category $D^b_X(\cL^{-1})$, and we claim that the global statement follows (but since this discussion is only motivational we skip the argument). 

Since $\cE$ is just the sheafy analogue of the extension algebra for the `bimodule' $\cL^{-1}[-1]$, it's clear from this description 
that this is the same construction as in Example \ref{eg.algebraic}, up to our choice of how much we shifted the bimodule.

In particular we see how to remove the unwanted shift in Example \ref{eg.geometric}. We replace $\cL^{-1}$ by the `shifted line-bundle' $\cL^{-1}[2]$, which is a sheaf of graded polynomial algebras over $X$, and consider the derived category $\D_X(\cL^{-1}[2])$. This is locally generated by the push-foward of $\cO_X$, and by Koszul duality is equivalent to the derived category of the sheaf of graded algebras:
$$  \cO_X \oplus \left( \cL^{-1}[1]\right)$$
Then the twist around the functor $i_*: \D_X(\cL^{-1}[2]) \to \D(X)$ is precisely $\otimes \cL$.

\subsection{The general case}\label{sect.generalcase}

Now we abstract this construction. Let $\cD$ be any triangulated category, and let
$$\Phi : \cD \to \cD$$
be any autoequivalence. We will show that $\Phi$ can be formulated as a spherical twist. 

Firstly we need to construct a category playing the role of $\D(E)$ or $\D_X(\cL^{-1})$, but this is not difficult. We consider a category $\cD_\Phi$ which has the same objects as $\cD$, and for each object $x\in \cD$ we denote by $j^*x$ the corresponding object in $\cD_\Phi$. Then we define the  morphisms in $\cD_\Phi$  by

$$\Hom_{\cD_\Phi}\!\!\big(j^*x,j^*y\big) = \Hom_{\cD}\!\big(x,y\big)\oplus \Hom_\cD\!\big(x, \Phi^{-1} y[1]\big)$$
with the composition law:
\beq{eq.composition}(f, f')\circ (g, g') = \big(f\circ g,\; f'\circ g + \Phi^{-1}(f)\circ g' \big) \eeq

\begin{rem} This is really just Example \ref{eg.algebraic} again, done with categories instead of algebras. View the category $\cD$ as a bimodule over itself, with the standard right action, but with the left action twisted by $\Phi^{-1}[1]$. Then $\cD_\Phi$ is exactly the extension category of $\cD$ by this bimodule.
\end{rem}

\begin{rem}\label{rem.foundations2}As stated this category is not triangulated, but we can correct this with the following steps. Firstly, we must replace $\cD$ by a dg-enhancement $\cD_{dg}$, and assume that we can lift $\Phi$ to a quasi-functor. We need to assume that $\Phi$ and $\Phi^{-1}$ are strict dg-functors, but we can achieve this by taking a cofibrant model for $\cD_{dg}$ \cite{Tabuada}. Then exactly the same construction produces a dg-category $(\cD_{dg})_\Phi$ whose homotopy category is $\cD_\Phi$ (we need $\Phi^{-1}$ to be strict or the composition law will not be strictly associative). Finally, we can let $\overline{(\cD_{dg})_\Phi}$ be the closure of $(\cD_{dg})_\Phi$ under mapping cones, this can be constructed as the category of perfect complexes over $(\cD_{dg})_\Phi$.  Then the homotopy category of $\overline{(\cD_{dg})_\Phi}$ is our required triangulated category.
\end{rem}

\begin{rem}Suppose, as in Example \ref{eg.geometric}, that $\cD=\D(X)$ and $\Phi = \otimes \cL[k]$ for a line bundle $\cL$ on a scheme $X$ and some $k\in \Z$. Then the category $\cD_\Phi$ is exactly $\D_X(\cL^{-1}[2-k])$ (after applying the corrections of the previous remark).   If instead $\cD=\D(A)$ and $\Phi=\otimes B[k]$ as in Example \ref{eg.algebraic}, then $\cD_\Phi$ is the subcategory of $D^b(A\oplus B^\vee[1-k])$  generated by the image of the functor $j^*$.
\end{rem}

We denote by $j^*$ the obvious faithful functor:
$$j^*: \cD \to \cD_\Phi $$
This has a right adjoint
$$j_*:\cD_\Phi\to \cD$$
defined on the generating objects by:
$$j_*: \,j^*x \,\mapsto\, x \oplus \Phi^{-1} x[1] $$

For the morphisms between these objects, we have that
$$j_*: \Hom_{\cD_\Phi}(j^*x, j^*y) \to \Hom_\cD(j_*j^*x,\,j_*j^*y)$$
is the map which sends
$$(f, f') \in  \Hom_\cD\!\big(x, \, y\oplus\Phi^{-1} y[1]\big)$$
to:
\begin{equation}\label{eqn:j_*}
\begin{pmatrix} f & f' \\ 0 &\Phi^{-1} f \end{pmatrix} \in  \Hom_\cD\!\big( x \oplus \Phi^{-1} x[1],\, y\oplus \Phi^{-1} y [1]\big)
\end{equation}

Then the right adjoint to $j_*$ is the functor 
$$ j^! = j^*\circ \Phi[-1]$$
 and it's automatic that the twist around $j_*$ is the autoequivalence $\Phi$. Thus conditions (i) and (iii) from Definition \ref{defn:spherical} are satisfied. Furthermore, all these functors immediately lift to strict dg-functors, as required by Remarks \ref{rem.foundations1} and \ref{rem.foundations2}. So the functor $j_*$ is spherical and $\Phi$ is a spherical twist. We have proved:
\begin{thm}\label{thm:mainthm} Every autoequivalence of $\cD$ arises as a spherical twist.\end{thm}

\begin{rem} \label{rem.monad}
Recall that given any monad $M: \cD \to \cD$, there are various ways to formally construct an adjunction giving rise to $M$. One way is to use the Kleisli category of `free $M$-algebras' \cite[p. 147]{Maclane}; here we construct a new category $\cD_M$ with the same objects as $\cD$, but with the morphisms between $x$ and $y$ replaced by $\Hom_\cD(x, My)$. The monad structure allows us to compose these morphisms, and there is an evident adjunction between $\cD$ and $\cD_M$ that recovers $M$. Our category $\cD_{\Phi}$ is simply the Kleisli category for the monad  $\id_\cD\oplus \Phi^{-1}[1]$.
\end{rem}

If $\Phi$ is originally given to us as a twist around some spherical functor 
$$F:\cC \to \cD$$
 then the above construction certainly does not recover $\cC$ and $F$. For example if $\Phi$ is a twist around a spherical object then $\cC=\D(pt)$, whereas $\cD_\Phi$ is very much larger. In the next section we discuss how to address this deficit.

\section{Autoequivalences with a section}\label{sect.withsection}

\subsection{More examples}\label{sect.moreexamples}

Let us consider some further examples, closely related to our examples from Section \ref{sect.2examples}.

\begin{eg}\label{eg.algebraic2}
Let $A$ be an algebra, and let $z\in A$ be a central element which is not a zero divisor. Then $z$ generates a two-sided ideal $(z)\subset A$, and the quotient $A/(z)$ is an algebra. We have an algebra map $\hat{\jmath} : A \to A/(z)$, and a functor
$$\hat{\jmath}_*: \D(A/(z)) \to \D(A)$$
with left and right adjoints 
$$\hat{\jmath}^*=A/(z)\otimes -\aand \hat{\jmath}^! = \Hom_A(A/(z), -)$$
 (both these functors are derived).  Using the Koszul resolution 
\beq{eq.koszulofAmodz}A \stackrel{z}{\To} A \To A/(z) \eeq
it's immediate that $\hat{\jmath}^! = \hat{\jmath}^*[-1]$, and that the functor $\hat{\jmath}_* \hat{\jmath}^*$ agrees with tensoring by the dg-$A$-bimodule $[A\stackrel{z}{\to} A]$. It follows easily that 
 $\hat{\jmath}^*$ is spherical and that the twist around it is the identity functor.

Since the identity functor can be described as tensoring with the diagonal bimodule $A$, we know from Example \ref{eg.algebraic} that we could instead have produced this autoequivalence using the graded algebra $E = A\oplus (A[1])$, and the spherical functor:
$$j_*: \D(E)\to \D(A)$$
 The relationship between these two constructions is not hard to see - the Koszul resolution  \eqref{eq.koszulofAmodz} in fact produces a dga, which is quasi-isomorphic to $A/(z)$. This dga is a deformation of $E$, and the map $j: A \to E$ obviously deforms to a map of dga's, lifting $\hat{\jmath}$.
\end{eg}

\begin{eg}\label{eg.algebraic3}
Now let $B$ be an arbitrary invertible $A$-bimodule, and let $z\in B$ be a central element, \emph{i.e.}~a bimodule map:
$$z: A \to B $$
 We saw in Example \ref{eg.algebraic} that the autoequivalence $\otimes B$ is a spherical twist around a spherical functor $$j_*: \D(E) \to \D(A)$$ where $E$ is the graded extension algebra:
$$E = A\oplus \left(B^\vee[1]\right) $$
Using $z$, we can deform $E$ to a dg-$A$-bimodule
$$E^z = \left[ B^\vee \stackrel{z}{\To} A \right]$$
(by which we mean keep the $A$-bimodule structure the same, and simply add this differential).  In the previous example, where we had $B=A$, this object was automatically a dga. However in general it is not, because the differential may not be a derivation. For it to be a derivation, we require additionally that 
$$z(\beta_1)\beta_2  = \beta_1 z(\beta_2), \quad\quad\forall \beta_1, \beta_2\in B^\vee $$
(since $\beta_1\beta_2=0$ by definition). Equivalently, we require the equality
\beq{eq.conditiononz} b\otimes z = z\otimes b \quad \in B\otimes_A B \eeq
for any $b\in B$. Let us assume this condition, so $E^z$ is a dga. The map $j$ obviously deforms to a map of dga's $j: A \to A^z$, and it's easy to calculate that the functor
$$j^*: \D(E^z) \to \D(A)$$
is still spherical, and the associated spherical twist is still $\otimes B$. Hence given such a $z$, we have a second way to produce this autoequivalence as a spherical twist.
\end{eg}

For completeness, we also describe the geometric analogue of the above construction.

\begin{eg}\label{eg.geometric2}

Let $X$ be a scheme and  let $\cL$ be a line-bundle on $X$. Suppose we have a non-zero section $\sigma\in H^0(\cL)$, whose zero-locus is a divisor: 
$$j:Y\into X$$
The push-forward functor
$$j_*: \D(Y) \to \D(X)$$
 is a classic example of a spherical functor (see \cite{Add}). Using the Koszul resolution
\beq{eq.KoszulofOY}\cL^{-1}\stackrel{\sigma}{\To} \cO_X \To j_*\cO_Y\eeq
it is easy to compute that the associated spherical twist is precisely $\otimes \cL$. 

We saw in Example \ref{eg.geometric} that we can also produce this autoequivalence from the spherical functor:
$$\iota^*: \D_X(\cL^{-1}[2]) \to \D(X)$$
To relate these two constructions, we recall that $\D_X(\cL^{-1}[2])$ can be viewed (via Koszul duality) as the derived category of the sheaf of graded algebras $\cO_X \oplus (\cL^{-1}[1])$. The Koszul resolution \eqref{eq.KoszulofOY} can be viewed as a sheaf of dga's over $X$, it is the `derived zero locus' of $z$. This sheaf of dga's is equivalent to $\cO_Y$, viewing the latter as a sheaf of algebras over $X$, and it is a deformation of $\cO_X \oplus (\cL^{-1}[1])$ (which is the derived zero locus of the zero section).

Note that in this example (just as in Example \ref{eg.algebraic2}) the analogue of condition \eqref{eq.conditiononz} holds automatically.
\end{eg}

\subsection{Wishful thinking}

As the examples of the previous section illustrate, we should consider a situation in which we have an autoequivalence $\Phi$ of some dg-category $\cD$, and also have a \emph{section} of $\Phi$, \emph{i.e.} a natural transformation:
$$\sigma: \id_\cD \to \Phi$$
Given such a $\sigma$ we might hope to construct a second spherical functor, whose twist is also $\Phi$, as a deformation (in some sense) of the functor $j_*: \cD_\Phi\to \cD$ that we constructed in Section \ref{sect.general}.

However there is a problem here, even at the purely formal level. If we set $\cD=\D(A)$ and $\Phi=\otimes B$ as in Example \ref{eg.algebraic3}, then a natural transformation from $\id_\cD$ to $\Phi$ is just a bimodule map $z: A \to B$. We saw that we needed the extra condition \eqref{eq.conditiononz}, which states that $z$ is a derivation. 

Let us translate this condition to our general situation. For any object $x\in \cD$, the component of $\sigma$ at $x$ is a morphism $\sigma_x: x \to \Phi(x)$. Applying $\Phi$, we get a morphism $\Phi(\sigma_x): \Phi(x) \to \Phi^2(x)$. But there is another evident morphism between these two objects, namely the component $\sigma_{\Phi(x)}$ of $\sigma$ at the object $\Phi(x)$. The analogue of condition \eqref{eq.conditiononz} is the statement that
\beq{eq.conditiononsigma} \Phi(\sigma_x) = \sigma_{\Phi(x)} \quad\in \Hom_\cD(\Phi(x), \Phi^2(x)) \eeq
for any $x\in \cD$.\footnote{It may be tempting to think that this condition holds automatically for any natural tranformation! Of course it does not, as Example \ref{eg.algebraic3} shows.} To put it another way: recall  from Remark \ref{rem.monad} that the endofunctor
$$\id_\cD \oplus \Phi^{-1}[1] $$
is a monad, \emph{i.e.} an algebra object in the monoidal category of endofunctors of $\cD$. Condition \eqref{eq.conditiononsigma} is the statement that the natural transformation
$$ \sigma: \Phi^{-1} \to \id_\cD $$
is a derivation. Consequently the cone on $\sigma$ defines a dg-algebra in endofunctors of $\cD$, which is a deformation of $\id_\cD \oplus \Phi^{-1}[1]$. 

At a formal level it is now clear what to do, we just apply the Kleisli construction to this new monad. That is, we define a category $\cD_\Phi^\sigma$ which has the same objects as $\cD$, and has morphism spaces
$$\Hom_{\cD_\Phi^\sigma}(j^*x, j^*y)  =\Hom_\cD\!\left( x,\;   \left[   \Phi^{-1} y   \stackrel{\sigma}{\To} y \right] \right) $$
where as before we use $j^*x$ to denote the object of $\cD_\Phi^\sigma$ corresponding to $x\in \cD$. Composition is still defined by the rule \eqref{eq.composition}, and the condition \eqref{eq.conditiononsigma} guarantees that the differential is a derivation. Evidently this category $\cD_\Phi^\sigma$ is a deformation of $\cD_\Phi$, in some sense.

  It is easy to check that the functors $j_*$ and $j^*$ extend to this deformed category and remain adjoint, that $j_*$ is still spherical, and that the twist around $j_*$ is still the autoequivalence $\Phi$ (these calculations are exactly the same as those appearing in the examples of Section \ref{sect.moreexamples}).

\begin{rem}\label{rem.strictness}
Unfortunately, for all of this to work as stated one must assume not only that $\Phi^{-1}$ is a strict dg-functor, but that $\sigma$ is a strict natural transformation, and that the condition \eqref{eq.conditiononsigma} holds strictly. We don't know any way to strictify all three requirements simultaneously. 
 
Alternatively, if one was sufficiently versed in the foundations of dg-categories or stable $\infty$-categories, one could proceed as follows. Instead of asking for these conditions to hold strictly, we just ask for them to hold up to homotopy, and to be given the data of all necessary coherencies. This means that we are still asking for an algebra structure on the object $\id_\cD\oplus \Phi^{-1}[1]$, but it is now something like an $A_\infty$-algebra, or an algebra in a monoidal $\infty$-category. Then we could  take the Kleisli category of this monad, in the appropriate homotopical sense. However, although it seems likely that one could make this definition precise, it is not clear how you would ever produce a useful example of this structure.
\end{rem}

Let us pretend we haven't read the previous remark, and proceed. Notice that if $x\in\cD$ is an object such that $\sigma:x\to \Phi x$ is a homotopy-equivalence, then $j^*x$ is contractible in $\cD_\Phi^\sigma$ (in Example \ref{eg.geometric2} this is just the statement that any object supported away from the zero locus of $\sigma$ goes to zero in $\D(Y)$). This means that $\cD_\Phi^\sigma$ is typically `smaller' than $\cD$, whereas the undeformed category $\cD_\Phi$ was `larger'. So our deformed spherical functor
$$j_*: \cD_\Phi^\sigma \to \cD$$
seems like a more efficient way to construct $\Phi$ as a spherical twist. 
\pgap

Now we explain how this construction solves (formally!) the problem raised at the end of Section \ref{sect.general}, of reconstructing a spherical functor from the spherical twist. Suppose we are given a spherical functor
$$F: \cC \to \cD$$
with right adjoint $R$ and counit $\epsilon: FR \to \id_\cD$. Let $\Phi$ be the spherical twist:
$$\Phi = \left[ FR \stackrel{\epsilon}{\To} \id_\cD \right] $$ 
We let $\sigma$ be the obvious natural transformation:
$$\sigma: \id_\cD \to \Phi$$
We claim that condition \eqref{eq.conditiononsigma} holds for this $\sigma$, so we have a category $\cD_\Phi^\sigma$ and spherical functor $j_*: \cD_\Phi^\sigma \to \cD$.

\begin{rem}\label{rem.strictness2}The condition certainly does not hold strictly, but we claim that it does hold up to homotopy. The argument for this claim is slightly fiddly, and given Remark \ref{rem.strictness} it doesn't seem worthwhile wasting the reader's time with it.
\end{rem}

Now $\Phi^{-1}$ is given by the mapping cone
$$\Phi^{-1} = \left[ \id_\cD  \stackrel{\eta}{\To} FL\right][-1]$$
where $L$ is the left adjoint to $F$, and $\eta$ is the unit of the adjunction. Recall that $\cD_\Phi^\sigma$ arises from considering the monad:
$$ \left[ \Phi^{-1} \stackrel{\sigma}{\To} \id_\cD \right]$$
This monad is obviously quasi-isomorphic to the endofunctor $FL$ (morally, this is the reason that $\sigma$ satisfies \eqref{eq.conditiononsigma}). In terms of categories, this is the observation that the morphisms in $\cD_\Phi^\sigma$ are given by
\begin{align}  \nonumber
\Hom_{\cD_\Phi^\sigma}(j^*x, j^*y) & := \Hom_\cD\!\left( x,\;   \left[   \Phi^{-1} y  \stackrel{\sigma}{\To} y \right] \right) \\
\label{eq.problem} &\simeq  \Hom_\cD\!\left( x,  FLy \right)\; \simeq\; \Hom_\cC\!\left( Lx,  Ly\right)  \end{align}
for each pair of objects $x,y\in \cD$. So the category $\cD_\Phi^\sigma$ is equivalent to the essential image of the functor $L$ inside the category $\cC$.

Furthermore, since the functor $j^*: \cD \to \cD_\Phi^\sigma$ is just $L: \cD \to \Im(L)$, passing to right adjoints shows that our spherical functor $j_*: \cD_\Phi^\sigma\to \cD$ must be the same as the given spherical functor $F$ (restricted to $\Im(L)$). So we have shown:

\begin{'prop'}\label{prop.reconstruct} Let $F: \cC \to \cD$ be a spherical functor whose left adjoint $L$ is essentially-surjective. Let $\Phi$ be the associated spherical twist and $\sigma:\id_\cD \to \Phi$ be the obvious natural transformation. Then we have an equivalence
$$\cC\, \cong\, \cD_\Phi^\sigma$$
under which the functors $F$ and $j_*$ correspond.
\end{'prop'}

Of course we haven't really proved this result, since we're ignoring the problems raised in Remarks \ref{rem.strictness} and \ref{rem.strictness2}.

Notice that without the essential-surjectivity condition we cannot hope to reconstruct $\cC$ from $\Phi$. For example if $\cC'$ is any category at all, and $\pi: \cC\oplus \cC'\to \cC$ is the projection functor, then
$$ F \pi: \cC\oplus \cC' \to \cD$$
is also a spherical functor whose effect on $\cD$ is indistinguishable from that of $F$.
\pgap

\noindent\textbf{Erratum.} Unfortunately the argument above is not just technically incomplete, it is actually wrong. We have left it in place since it appears in the published version of this paper, but there is a very simple counter-example which we learnt from \cite{Merlin}. 

Let $\cC$ be the derived category of the semi-simple algebra $\C\oplus \C$ and let $\cD=D^b(\C)$. We have a homorphism $i: \C \to \C\oplus \C$  and a pair of bi-adjoint functors $F=i_*: \cC \to \cD$ and $L = i^*: \cD\to \cC$. It's easy to check that $F$ is spherical and the twist around it is the shift $\Phi = [1]$. But the category $\cD_\Phi^\sigma$ is not $\cC$, it's the derived category of the nilpotent algebra $\C[\epsilon]/\epsilon^2$.

The flaw in our argument is that we do not have a functor between $\cC$ and $\cD_\Phi^\sigma$, so \eqref{eq.problem} is only an equality of vector spaces and need not respect the composition of morphisms. This is very clear in the counterexample since:
$$\Hom_{\cD_\Phi^\sigma}(j^*\C, j^*\C) = \C[\epsilon]/\epsilon^2 \;\neq\; \C\oplus \C = \Hom_{\cC}(L\C, L\C) $$

The conclusion is that $\cC$ probably is a deformation of $\cD_\Phi^\sigma$, but in general one must deform the multiplication on the monad as well as the differential.

\section{$\P$-twists}\label{sect.Ptwists}

In this section (which is logically independent of the rest of the paper) we'll consider the $\P$-twist autoequivalences defined by Huybrechts and Thomas \cite{HuyTho}. As we'll see, these can be reformulated as spherical twists very explicitly.

For this section we restrict to the case that our target category $\cD = \D(X)$
 is the derived category of a smooth projective variety $X$ over $\C$. We could use other $\cD$ if we assume suitable finiteness conditions, Serre duality, etc. 

\begin{defn}\cite{HuyTho} An object $P\in \D(X)$ is called a $\P^{n}$-object if we have $\omega_X\otimes P \cong P$ and 
$$\Ext^\bullet(P, P ) \cong \C[h] / (h^{n+1}) $$
 as a graded ring, where the generator $h$ has degree 2.
\end{defn}
By Serre duality, a $\P^n$-object can only exist when $\dim X = 2n$.  Huybrechts and Thomas used $\P^n$-objects to construct a new kind of autoequivalence of a derived category, a `$\P$-twist', whose kernel  is the double cone
\begin{center}\begin{tikzcd}[column sep = huge]
P^\vee\boxtimes P [-2] \arrow{r}{h\otimes 1 - 1\otimes h} &P^\vee\boxtimes P \arrow{r}& \cO_{\Delta}
\end{tikzcd}\end{center}
in $\D(X\times X)$. We will show that this definition can be naturally interpreted as a spherical twist. 

Let $A$ be the graded algebra
$$ A = \C[h]$$
with the generator $h$ having degree $2$. We view $A$ as a dga with zero differential. Suppose we have a $\P^n$-object $P\in \D(X)$. The object $P$ carries an action of $A$, and this defines a functor
$$F = P\otimes_A - : \;\D(A) \to \D(X)$$
under which $A$ maps to $P$. Note that any object in $D^b(A)$ has a bounded free resolution,  so this functor does land in $\D(X)$ and not some larger category.
\begin{prop} Assume that $\Hom(P,P)$ is formal as a dg-module over $A$. Then the functor $F$ is spherical, and the spherical twist around $F$ is the $\P$-twist around $P$.
\end{prop}
 We don't believe that the formality hypothesis is essential, more likely it is just an artefact of this paper's persistently slapdash attitude to foundations.
\begin{proof}
The right adjoint of $F$ is 
$$R= \Hom(P, -): \D(X) \to \D(A)$$
For any object $\cE\in \D(X)$ the homology of the complex $R(\cE)$ is finite-dimensional over $\C$,  so consists of torsion $A$-modules. This means we can apply Serre duality for $A$ and $X$ to deduce\footnote{Note that $A$ behaves like a non-compact Calabi-Yau of dimension $-1$, in that the diagonal bimodule satisfies $\Hom_{A-A}(A, A) = A[1]$.} that the left adjoint of $F$ is:
$$L =  \Hom(\omega_X\otimes P, -)[2n+1] $$
Since $\omega_X\otimes P \cong P$, this means that $L=R[2n+1]$. The cotwist associated to $F$ is the functor:
$$ C = \big[1 \to RF \big][-1]$$
For the module $A\in D^b(A)$, we have
\begin{equation}\label{eq.formal}RF(A) = \Hom(P,P) \simeq A/(h^{n+1}) \end{equation}
(using our formality hypothesis). From the exact sequence
\begin{center}\begin{tikzcd}[column sep = 30pt]
 A[-2n-2] \arrow{r}{h^{n+1}} &A \arrow{r} & A/(h^{n+1}) 
\end{tikzcd}\end{center}
we conclude that $C(A)=A[-2n-2]$. Furthermore, if we take the endomorphism $h\in \Hom_{\D(A)}(A,A)$ then it's clear that $C(h) = h$. Since $A$ generates $D^b(A)$, this proves that $C$ is just the shift $[-2n-2]$.  Thus conditions (ii) and (iv) of Definition \ref{defn:spherical} are satisfied and $F$ is spherical.

Now we examine the spherical twist around $F$, which is the cone $[FR \to 1]$. This sends an object $\cE\in \D(X)$ to:
$$\big[ P\otimes^L_A \Hom(P, \cE) \To \cE\big] $$ 
Using the free resolution of the diagonal $A$-bimodule
$$\begin{tikzcd}[column sep = huge]
 A\otimes A [-2] \arrow{r}{h\otimes 1 - 1\otimes h} & A\otimes A 
\end{tikzcd}$$
we have that  $  P\otimes^L_A \Hom(P, \cE) $ is equal to the cone on
$$\begin{tikzcd}[column sep = huge]
 P\otimes\Hom(P,\cE)[-2]   \arrow{r}{h\otimes 1 - 1\otimes h} & P\otimes \Hom(P,\cE) 
\end{tikzcd}$$
and we recover the formula for the $\P$-twist.
\end{proof}

\begin{rem} As before, we don't really need the whole category $\D(A)$ but only the essential image of the functor $R$, which is the category of torsion modules.  By Koszul duality, this subcategory is equivalent to the category of perfect complexes over the algebra
$$B=\C[\epsilon]/(\epsilon^2)$$
where $\epsilon$ has degree $-1$. Then $F$ sends the object $B\in \Perf(B)$ to the object:
\beq{eq.S}S = \big[ P[-1] \stackrel{h}{\To} P \big] \quad \in \D(X) \eeq
This object $S$ is similar to a `fat spherical object' in the sense of Toda, but over an `odd' analogue of the dual numbers.
\end{rem}

\begin{rem}\label{rem.Pfunctors}
In \cite{Add}, Addington introduces a generalization of $\P^n$-objects, to $\P^n$-\emph{functors}. These involve a source category $\cC$, an autoequivalence $H$ of $\cC$, and a functor
$$F : \cC \to \cD$$
with right and left adjoints $R$ and $L$. There are three axioms, of which the most fundamental is that:
$$RF= \id_\cC \oplus H \oplus H^2 \oplus... \oplus H^n $$
If $\cC= \D(pt)$ and $H=[-2]$ then the definition reduces to a choice of $\P^n$-object.

It seems that one can extend the construction given above for $\P^n$-objects to $\P^n$-functors, in the following way. Set $H' = H[2]$, and consider the category $\cC_{H'}$ as constructed in Section \ref{sect.general}.  In the case $\cC = \D(pt)$ and $H=[-2]$, this category $\cC_{H'}$ is exactly $\Perf(B)$, where $B=\C[\epsilon]/(\epsilon^2)$ as in the previous remark. We have a functor $j_*: \cC_{H'}\to \cC$, and to extend what we did for $\P$-objects we should consider the functor:
\al{ FHj_*[1] : \cC_{H'} & \To \cD \\ j^*x & \mapsto FHx[1] \oplus Fx }
Because $F$ is a $\P$-functor, we have in particular a natural transformation from $H$ to $FR$, which induces a natural transformation $h:FH\to F$ by adjunction. We can deform $FHj_*[1]$ to get a functor $\Sigma: \cC_{H'} \to \cD$, with:
$$\Sigma: j^*x \to \left[ FHx \stackrel{h}{\To} Fx \right]$$
In the case of $\P$-objects, this functor is the object $S$ \eqref{eq.S}.

Presumably $\Sigma$ is spherical. It is not so easy to check this assertion, because the full definition of a $\P$-functor is a little involved. However, this construction could be taken as an alternative (and perhaps more general) definition of a $\P$-functor. This approach is being pursued by Anno and Logvinenko \cite{AnnLogIP}.
\end{rem}

\vspace{1cm}

Ed Segal

University College London

\emph{e.segal@ucl.ac.uk}

\end{document}